\documentclass[lettersize,journal]{IEEEtran}
\usepackage{amsmath,amsfonts}
\usepackage{algorithmic}
\usepackage{algorithm}
\usepackage{array}
\usepackage[caption=false,font=normalsize,labelfont=sf,textfont=sf]{subfig}
\usepackage{textcomp}
\usepackage{stfloats}
\usepackage{url}
\usepackage{verbatim}
\usepackage{graphicx}
\usepackage{cite}
\usepackage{booktabs}
\usepackage{graphicx}
\usepackage{caption}
\usepackage{multirow}
\usepackage{float}

\begin{document}

\title{Efficient Link Prediction in Continuous-Time\\ Dynamic Networks using Optimal Transmission and Metropolis Hastings Sampling}

\author{Ruizhi Zhang, Wei Wei, Qiming Yang, Zhenyu Shi, Xiangnan Feng, Zhiming Zheng
\thanks{Ruizhi Zhang, Qiming Yang and Zhenyu Shi are with the School of Mathematical Sciences, Beihang University, 100191, Beijing, China.}
\thanks{Wei Wei is with the School of Mathematical Sciences, Beihang University,
	100191, Beijing, China, Key Laboratory of Mathematics Informatics
	Behavioral Semantics, Ministry of Education, Beijing, 100191, China,
	Institute of Artificial Intelligence, Beihang University, 100083, Beijing,
	China, and Zhongguancun Laboratory, Beijing, 100094, China.}
\thanks{Xiangnan Feng is with the Complexity Science Hub, Vienna, 1080, Austria.}
\thanks{Zhiming Zheng is with the Institute of Artificial Intelligence, Beihang University, 100083, Beijing,China.}
\thanks{Corresponding author: Wei Wei (weiw@buaa.edu.cn).}

}
\markboth{Journal of \LaTeX\ Class Files,~Vol.~14, No.~8, August~2021}%
{Shell \MakeLowercase{\textit{et al.}}: A Sample Article Using IEEEtran.cls for IEEE Journals}



\maketitle
\begin{abstract}
Efficient link prediction in continuous-time dynamic networks is a challenging problem that has attracted much research attention in recent years. A widely used approach to dynamic network link prediction is to extract the local structure of the target link through temporal random walk on the network and learn node features using a coding model. However, this approach often assumes that candidate temporal neighbors follow some certain types of distributions, which may be inappropriate for real-world networks, thereby incurring information loss. To address this limitation, we propose a framework in \underline{c}ontinuous-time dynamic networks based on \underline{O}ptimal Transmission (OT) and \underline{M}etropolis Hastings (MH) sampling (COM). Specifically, we use optimal transmission theory to calculate the Wasserstein distance between the current node and the time-valid candidate neighbors to minimize information loss in node information propagation. Additionally, we employ the MH algorithm to obtain higher-order structural relationships in the vicinity of the target link, as it is a Markov Chain Monte Carlo method and can flexibly simulate target distributions with complex patterns. We demonstrate the effectiveness of our proposed method through experiments on eight datasets from different fields.
\end{abstract}


\section{Introduction}
\IEEEPARstart{N}{etwork} data is ubiquitous across numerous fields, including social networks, academic citation and collaboration networks \cite{social_network, citation_network, collaboration_network}. The majority of links in these networks are continuously evolving over time, making it crucial to develop an effective continuous time dynamic network link prediction model. Such a model can reveal the underlying evolution process of the network, facilitating a better understanding of the structure, dynamics, and properties of the network. Furthermore, it can be used for various practical applications, including recommendation systems, traffic-flow prediction, and fraud detection \cite{recommender_system,traffic_flow, fraud_detecting}.

Heuristic link prediction algorithms have gained popularity due to their simplicity and ease of implementation \cite{CN, AA}, which are based on the assumption that nodes with similar properties are more likely to be connected \cite{similarity_assumption}. While these methods can provide useful insights in some cases, they have their limitations. Many of the indicators used in these algorithms are designed for specific tasks and may not be generalizable to real-world networks. Moreover, these methods may not always perform well on large-scale networks or networks with complex topologies.

To address these limitations, researchers have proposed more sophisticated approaches such as the random walk-based method, which uses biased random walks to learn the low-dimensional representation of nodes \cite{Node2vec, deepwalk}. Another promising approach is deep learning-based \cite{gcn, gan}, which aggregates neighbor information to extract important features automatically. These methods can perform well in some scenarios, but they often lack explanatory power and do not fully consider the temporal characteristics of the network.

Some researchers have addressed temporal information in networks by partitioning them into discrete time snapshots \cite{network_snapshot1, network_snapshot2}. However, this approach can lead to a loss of valuable temporal information and the selection of appropriate time partitions can be challenging. To overcome this, CTDNE \cite{CTDNE} treats the network as a continuous-time network, which helps to preserve temporal dynamics. The algorithm uses a continuous-time random walk to generate node sequences and learns node embeddings based on these sequences. 

Although CTDNE has demonstrated promising results on some networks, its effectiveness relies on the sampling strategy used in the temporal dimensions, assuming unbiased or biased distributions in the network. Despite that, real network data may not follow these distribution, resulting in information loss during sampling. To tackle this challenge, we propose a novel sampling approach that leverages optimal transmission, effectively utilizing a time-valid space for neighbor sampling \cite{Optimal_Transmission}. We calculate the Wasserstein \cite{Wasserstein_Distance} distance between the current node and candidate nodes using nodes probability distribution, ensuring selection of the node with the smallest distance and similar information propagation patterns for the next sampling \cite{Wasserstein_Distance}. Efficiently incorporating temporal and spatial information is vital for enhancing link prediction tasks. The plain random walk sampling method may yield irrelevant nodes, leading to an imprecise representation of the network's topology. To address this issue, we employ the Metropolis-Hastings (MH) algorithm to sample the local neighbor space structure of the target link \cite{MH}, as it can effectively handle distributions with complex patterns, including those with multiple modes or non-standard shapes. This makes it suitable for capturing the complexity and variability present in real-world data. Therefore, our proposed link prediction framework comprehensively captures information from both the temporal and spatial dimensions. 
Our contributions can be summarized as follows:
\begin{enumerate}
	\item{We propose a novel sampling strategy that includes the temporal dimension based on the optimal transmission theory, avoiding information loss resulting from solely considering the time validity of the next node. Additionally, we utilize MH sampling to acquire the node's high-order neighborhood structure.} 
	
	\item{We propose a new continuous time dynamic network link prediction scalable framework COM, which simultaneously exploits network temporal information and higher-order topological information to improve prediction results.}

	\item{The effectiveness of the proposed method was compared on eight real networks, and the results indicate its superiority over several state-of-the-art methods.}
	
\end{enumerate}

The remainder of this article is organized as follows. Section II reviews some related work and problem statements are introduced in Section III. Section IV presents the proposed algorithm in detail. Section V evaluates the performance of the COM, and Section VI concludes the investigation.
\section{Related Work}
Dynamic network link prediction is a critical problem in network analysis, particularly in the era of big data where the amount of information generated by networks is increasing rapidly. Existing link prediction approaches can be broadly categorized into three main paradigms: path-based methods, embedding methods, and graph neural networks.

\noindent \textbf{Path-based Methods.} Path-based methods are one of the traditional approaches for link prediction in dynamic networks. They aim to capture the similarity between nodes based on their common paths or walks. Commonly used path-based methods include Katz index \cite{katz}, and Local Path index \cite{LPI}. Several studies have proposed variations and extensions of path-based methods, such as the relation strength similarity (RSS) \cite{RSS} and the PathSim algorithm \cite{pathsim}. RSS captures potential relationships in real-world network structures, allowing users to specify relation strength and adjust the discovery range, while PathSim measures the similarity between two nodes by comparing the similarity of their paths. Path-based methods for link prediction offer high accuracy and interpretability, but suffer from computational complexity issues and the inability to capture dynamic network properties.

\noindent \textbf{Embedding Methods.} Link prediction using embedding methods has gained much attention in recent years due to their ability to capture both structural and semantic information in networks. These methods generate low-dimensional vector representations for nodes, which can be used for predicting links. One popular embedding method is the node2vec \cite{Node2vec} algorithm, which learns node representations by performing biased random walks on the network and then uses Skip-gram models to learn representations based on the sequences of nodes visited during the walks. Other methods, such as DeepWalk \cite{deepwalk} and LINE \cite{Line}, have also been proposed for link prediction based on node embeddings. Embedding-based methods for link prediction have the advantage of being computationally efficient and able to handle large-scale networks. However, they may struggle to capture complex network structures and dynamic network properties, leading to suboptimal performance in certain scenarios.

\noindent \textbf{Graph Neural Networks Methods.} Deep learning models can learn complex patterns and non-linear relationships from network data, and have achieved promising results in link prediction tasks. The most popular deep learning models for link prediction include graph convolutional networks (GCNs) \cite{gcn}, and attention-based models \cite{gan}. These models can be applied to various types of networks, such as social networks, biological networks, and knowledge graphs. TDGNN \cite{TDGNN} incorporates network temporal information into GNNs via a novel Temporal Aggregator to obtain target node representations by assigning neighbor nodes aggregation weights using an exponential distribution to bias different edges’ temporal information. GC-LSTM \cite{GC-LSTM} is a novel end-to-end model for dynamic network link prediction that combines a Graph Convolution Network (GCN) and LSTM to respectively capture local structural properties and temporal features of all snapshots of a dynamic network.
However, these methods may encounter challenges when dealing with large-scale networks and avoiding overfitting.

\section{Preliminaries AND Problem Formulation}
\noindent \textbf{Definition 3.1 (Continuous-Time Dynamic Networks).}
\textit{Continuous-time dynamic networks are networks in which the links between nodes can appear and disappear over time in a continuous fashion, rather than in discrete time steps. In mathematical terms, given a time-varying graph $\mathcal{G = (V, E_{T}, T)}$, where $\mathcal{V}$ is a set of nodes or vertices, and $\mathcal{E_{T}}$ is a set of edges or links that varies with time $\mathcal{T}$. $|\mathcal{V}|$ and $|\mathcal{E_{T}}|$ denote number of nodes and edges.} This paper concerns graphs that are directed and do not have weights assigned to the edges.

The edges in $\mathcal{E_{T}}$ can appear or disappear at any point in time, and are represented by a triplet $e^{t}_{i,j}=(v_{i}, v_{j}, t) \in \mathcal{E_{T}}$ , where $v_{i}$ and $v_{j}$ are the nodes that the edge connects. Nodes that are connected imply that they are neighbors to each other. We further define the temporal neighbors.\\

\noindent \textbf{Definition 3.2 (Temporal Neighborhood).}
\textit{For a given node $v \in \mathcal{V}$ and a timestamp $t \in \mathcal{T}$, the temporal neighborhood of $v$ is denoted as $\mathcal{N}_{t}(v) $. Formally,}
\begin{equation}
	\label{temporal_neighbor}
	\mathcal{N}_{t}(v) = \{(w,t^{\prime}) |e^{t^{\prime}}_{v,w}=(v,w,t^{\prime}) \in \mathcal{E_{T}} \land t^{\prime} > t\}.
\end{equation}
The size of the temporal neighborhood of a node $v$ at time $t$ is given by the cardinality of $\mathcal{N}_{t}(v)$, denoted as $|\mathcal{N}_{t}(v)|$. The temporal neighborhood provides information about the immediate network context of a given node, which can be useful in predicting its future connections.\\

\noindent \textbf{Definition 3.3 (Continuous-Time Dynamic Networks Link Prediction).}
\textit{Continuous-time dynamic network link prediction is a subfield of dynamic network link prediction that focuses on predicting the appearance or disappearance of links over a continuous time period. This can be formulated as a random event, where the occurrence of each link is a random event that happens at a specific time with a certain probability. The goal of link prediction is to predict the presence or absence of a link between any two nodes $v_{i}$ and $v_{j}$ in the future graph $\mathcal{G = (V, E_{T}, T +}$ $\Delta T)$, where $\Delta T$ is the time horizon.}

\noindent \textbf{Definition 3.4 (Temporal Network Embedding).}
\textit{The goal of temporal network embedding is to learn a function that maps each node in the network to a low-dimensional vector that captures its temporal dynamics and network topology. Formally, mapping $\phi : \mathcal{V} \to \mathbb{R}^{l}, \phi(u)=\emph{x}_{u}$, where $\emph{x}_{u}$ is called an embedding of node $u \in \mathcal{V}$.}

\section{Methodology}
\begin{figure*}[h!t]
	\centering
	\includegraphics[width=0.9\linewidth]{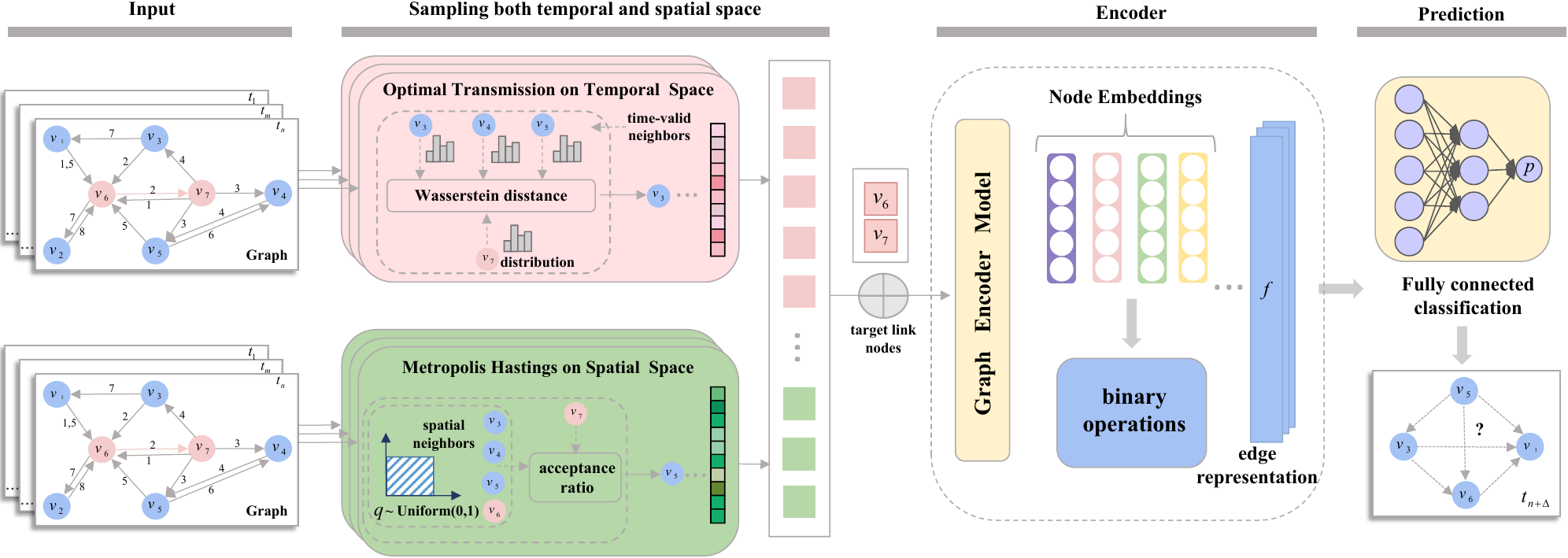}
	\caption{The overall framework of COM. COM first samples nodes in the graph from both the temporal and spatial dimensions. Specifically, the next node is sampled using the Wasserstein distance to ensure minimal information loss, and a Markov chain that follows the target distribution is obtained using the proposal distribution to obtain structural information. Node embeddings are then obtained using an encoding model, and edge representations are generated through binary operations. Finally, a fulled connected classifier is employed to score and predict potential edges.}
	\label{fig:framework}
\end{figure*}

In this section, we delve into the details of generating embeddings and their utilization for link prediction. Fig. \ref{fig:framework} illustrates an overview of the entire process. 
\subsection{Overview}
The key idea of COM is to establish a dynamic network link prediction framework that can capture continuously evolving spatiotemporal information on the graph. To achieve this, it is represented as an encoder-decoder model.
COM first extracts network patterns in both the temporal and spatial dimensions and encodes them as low dimensional vector. Then, the decoder combines the encoder with a simple MLP to make task-specific predictions. It can be represented by the following function.
\begin{equation}
	\mathcal{Z_{T}} = f(g(\mathcal{G(V, E_{T}, T)})).	
\end{equation}

The encoder $g$ in the equation above encodes the temporal graph to generate node representations based on optimal transmission and Metropolis Hastings sampling, which are then synthesized to obtain edge representations (the synthesis strategy will be analyzed in the experiment). The decoder $f$ predicts whether or not to connect edges based on vector decoding, and we use the simplest MLP for this task (other machine learning models could be used instead). In the equation, $\mathcal{Z_{T}} \in \{0,1\}$ represents the predicted value for the edge to be predicted, where the ground truth $\mathcal{Z_{T}}=1$ if there is an interaction between two nodes and 0 otherwise.

\subsection{Optimal Transmission and Temporal-biased Sampling}
\noindent \textbf{Optimal Transmission}
The French mathematician Monge first introduced the optimal transmission problem, which aimed to determine the most efficient way to transport sand piles from one location to another by considering various transmission methods with different associated costs. The objective of this problem is to identify the solution with the lowest transmission cost among all possible transmission schemes \cite{Optimal_Transmission}. 
Optimal transmission has been applied in various fields, including disease control and information dissemination. It has been shown to be an effective method for identifying the most influential nodes in a network and can lead to more efficient and effective information transmission.
Applying it to our tasks, optimal transmission refers to the process of transmitting information or data from one vertice to another in the most efficient and effective manner possible with minimal information loss. The optimal transmission problem is defined as minimizing the distance between the node and its neighbors,
\begin{equation}
	d(u, v) = \underset{v \in \mathcal{N}_{t}(u)}{min}<\mathcal{P}_{t}, \mathcal{G}>,		
\end{equation}
in which $\mathcal{N}_{t}(u)$ means the time-valid neighbors of $u$ at time $t$ and $\mathcal{P}_{t}$ represents the matrix composed of node probability distribution.

\noindent \textbf{Wasserstein Distance}
The optimal transmission problem's minimum value is represented by the Wasserstein distance, a distance measure in probability space that characterizes the minimum transmission cost from source node $u$ to target node $v$ \cite{Wasserstein_Distance}. It is a natural way to define distance between two probability measures in the measurement space, and provides a means to measure the difference between probability measures $u$ and $v$. The Wasserstein distance formal expression for the probability distributions of $u$ and $v$ at time $t$, assumed to be $p$ and $q$, is as follows:
\begin{equation}
	\mathcal{W}(p, q) = \underset{\gamma \sim \Pi (p, q)}{inf}E_{x,y \sim \gamma}[||x-y||].		
\end{equation}

The set of all possible joint distributions that combine distribution $p$ and $q$ is denoted by $\Pi(p, q)$. For each joint distribution $\gamma$ in this set, a pair of samples $x$ and $y$ can be sampled from it by drawing $(x, y) \sim \gamma$, and the distance $||x-y||$ between these samples can be calculated. The Wasserstein distance is defined as the infimum for the expected value of this distance over all possible joint distributions $\gamma$.

\noindent \textbf{Temporal-biased Sampling}
When it comes to node information propagation in networks, the sampling strategies of CTDNE rely on biased or unbiased distribution assumptions of nodes in a time-valid space, which can result in inaccuracies or loss of information. In contrast, the calculation of the Wasserstein distance based on the probability distribution of nodes offers a way to estimate the minimum loss of time information when sampling the next node. Specificly, the Wasserstein distance captures the similarity between the probability distributions of nodes in the network, enabling the identification of nodes with similar information propagation patterns and facilitating precise sampling of temporal neighbors.
More precisely, at a given time $t$, for a node $u$, we select the time-valid neighbors $\mathcal{N}_{t}(u)$ and calculate the Wasserstein distance based on the probability distribution $(\mathcal{P}_{t}(u), ..., \mathcal{P}_{t}(v_{i})(v_{i} \in \mathcal{N}_{t}(u)))$ of node $u$ and its neighbors at the current time $t$, and use it to determine the next sampling node of the random walk. We finally obtain the random walk sequence $\mathcal{W}_{t}(u) = (v_{1}, ..., v_{k}) (v_{1} \in \mathcal{N}_{t}(u),v_{i+1} \in \mathcal{N}_{t}(v_{i}))$ started from node $u$. In the time dimension, we limit the maximum step size of the random walk, without setting a specific walk length, as in the setting of CTDNE \cite{CTDNE}. The overall process of the proposed temporal-biased random walk strategy is given in Algorithm \ref{alg:temporal_sample}.
\begin{algorithm}
	\caption{Temporal-biased Sampling.}
	\begin{algorithmic}
		\STATE
		\STATE{\bf Input:}
		Input networks $\mathcal{G = (V, E_{T}, T)}$, max\_walk\_length $L$, start\_edge $e^{t}_{v,u}=(v,u,t)$, number of context windows $C$, context window size $ \omega $, number of walks for each node R, temporal context window count $ \beta  = R \times N  \times (L - \omega +1)  $.\\
		\STATE{\bf Output:} Temporal-biased random walk sequences $\textit{walks} \ \mathcal{W}_{t}(u)$.
		\STATE \hspace{0.5cm} Acquire $curr\_node = u$; $\textit{walks}$ = [ ]\\
		\STATE \hspace{0.5cm} $\textbf{while} \  \beta - C >0 \  \textbf{do}$
		\STATE \hspace{1cm} $\textbf{for} \ i = 1 : min\{L, C\}-1 \  \textbf{do} $
		\STATE \hspace{1.5cm} Compute time-valid neighbors $\mathcal{N}_{t}(curr\_node)$ with Eq. \ref{temporal_neighbor};\\
		\STATE \hspace{1.5cm} \textbf{if} $|\mathcal{N}_{t}(curr\_node)|>0$ \textbf{then}
		\STATE \hspace{2cm} $\textit{distance}$ = [ ]
		\STATE \hspace{2cm} \textbf{for} $v_{i} \in \mathcal{N}_{t}(curr\_node)$ \textbf{do}
		\STATE \hspace{2.5cm} $P_{t}(curr\_node) = p, P_{t}(v_{i}) = q$
		\STATE \hspace{2.5cm} $distance$.append($\textbf{Wasserstein}(p, q)$)
		\STATE \hspace{2cm} $index$ = argmin$(distance)$
		\STATE \hspace{2cm} $next\_node$ = $\mathcal{N}_{t}(curr\_node)[index]$
		\STATE \hspace{2cm} $curr\_node = next\_node$
		\STATE \hspace{2cm} $\textit{walks}$.append($curr\_node$)
		\STATE \hspace{1cm} $\textbf{if} \  |\textit{walks}| > \omega $
		\STATE \hspace{1.5cm} C = 	C + ($|\textit{walks}|$ - $\omega$ + 1) 
		\STATE \hspace{0.5cm} $\textbf{end while}$
	\end{algorithmic}
	\label{alg:temporal_sample}
\end{algorithm}

\begin{figure}[h!t]
	\centering
	\includegraphics[width=0.8\linewidth]{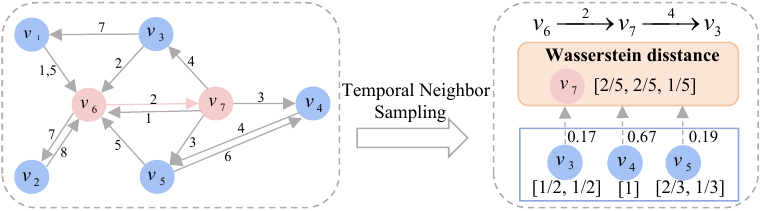}
	\caption{Schematic diagram of temporal space sampling based on optimal transmission. As shown in the figure, node $v_{7}$ is at time 2. Initially, a set of time-valid neighbors $\{v_{3},v_{4},v_{5}\}$ is selected to characterize the information distribution among nodes based on the degrees of their neighbors. Next, we calculate the Wasserstein distance between the current node $v_{7}$ and its neighboring nodes individually. The node $v_{3}$ with the smallest distance, indicating the least information propagation loss, is selected as the next sampling node.}
	\label{fig:temporal}
\end{figure}

\subsection{Metropolis Hastings and Spatial-biased Sampling}
\noindent \textbf{Metropolis Hastings Sampling}
Metropolis-Hastings (MH) Sampling is a Markov Chain Monte Carlo (MCMC) method, widely used for generating samples from a probability distribution that is difficult to sample directly \cite{MH_base}. Assuming $\pi (\cdot )$ is the target probability distribution defined on a sample space $\Omega$. The main idea behind MCMC is to construct a Markov chain on $\Omega$ that has $\pi (\cdot )$ as its stationary distribution. The MH algorithm is a specific MCMC method used to design a Markov chain with a desired stationary distribution.

To implement the MH algorithm, a conditional probability mass function $q(y|x)$ is first specified, where $x, y \in \Omega$. At each iteration $i$, a candidate sample $x^{c}$ is drawn from the proposal distribution $q(\cdot|x^{i-1})$, where $x^{i-1}$ is the sample selected in the previous iteration. The candidate sample is then accepted with a probability given by the MH acceptance ratio $p(x^{i-1}, x^{c})$, which depends on the ratio of the target distribution $\pi (\cdot )$ and the proposal distribution $q(\cdot|x^{i-1})$. The formula to compute the acceptance ratio is as follows,
\begin{equation}
	p(x^{i-1}, x^{c}) = min\{1, \frac{\pi(x^{c}) \cdot q(x^{i-1}|x^{c})}
		{\pi(x^{i-1}) \cdot q(x^{c}|x^{i-1})}\}.	
\end{equation}
If the candidate sample $x^{c}$ is accepted, then $x^{i}$ is updated to $x^{c}$. Otherwise, $x^{i}$ remains the same as the previous sample $x^{i-1}$.

\noindent \textbf{Spatial-biased Sampling}
Some graph tasks based on traditional spatial random walk methods perform deficiently, possibly due to their tendency to generate spatially unrelated samples, resulting in inaccurate representation of the underlying spatial structure of the network. Metropolis Hastings sampling overcomes this limitation by offering a more flexible proposal distribution, leading to improved sampling efficiency by reducing the number of rejected samples. Additionally, it can handle complex distributions, including those with multiple modes, that may be difficult for other random walk sampling methods. Clearly, we present a new edge sampler utilizing the Metropolis-Hastings (MH) sampling technique on graphs, where the target distribution is assumed to be the uniform distribution. To utilize the MH sampling method, we adopt the uniform distribution as the conditional probability mass function $q(\cdot|u)$, which is given by $q(\cdot|u)= 1/deg(u)$, where $deg(u)$ represents the degree of the node $u$. As a result of the symmetry of the uniform distribution, the acceptance ratio $p(x^{i-1}, x^{c}) $ for the $i$th candidate $x^{c}$ can be reduced to
\begin{equation}
	p(x^{i-1}, x^{c}) = min\{1, \frac{deg(x^{i-1})}{deg(x^{c})}\}.
\end{equation}
We choose the uniform distribution as the conditional probability mass function $q(\cdot|u)$ due to its computational efficiency in calculating acceptance rates. Furthermore, the use of uniform distribution for the edge sampler allows for convergence to any discrete target distribution\cite{MH}.

Our spatial-biased sampling strategy aims to generate a representative sampling set regardless of the distribution's complexity. By employing the MH algorithm, the resulting Markov chain is both aperiodic and recurrent, and converges to the target distribution. During the random walk based on the MH algorithm on the graph, given node $u$ at time $t$, the next node is sampled from all neighbors based on the MH algorithm. Unlike temporal-biased sampling, our sampling space includes all neighbors, not just the time-valid ones. Finally, we obtain the spatial-biased random walk sequence $\mathcal{W}_{s}(u) = (v_{1}, ..., v_{k}) (v_{1} \in \mathcal{N}(u),v_{i+1} \in \mathcal{N}(v_{i})) $ with a source node of $u$, where the length of the random walk is a hyperparameter. In our experiment, we set the step size to 8 based on several studies, which suggests that focusing on lower-order neighbours can provide more information for network structure \cite{lower-order_neighbor1, lower-order_neighbor2}. The overall procedure is described in Algorithm \ref{alg:spatial_sample}.
\begin{algorithm}
\caption{Spatial-biased Sampling.}
	\begin{algorithmic}
		\STATE 
		\STATE{\bf Input:}
		Input networks $\mathcal{G = (V, E_{T}, T)}$, max\_length $\textit{max\_length}$, start\_edge $e^{t}_{v,u}=(v,u,t)$.\\
		\STATE{\bf Output:} Spatial-biased random walk sequences $\textit{walks} \ \mathcal{W}_{s}(u)$.
		\STATE \hspace{0.5cm} Acquire $curr\_node = u$
		\STATE \hspace{0.5cm} $\textit{walks}$ = [ ]
		\STATE \hspace{0.5cm} $\textbf{while} \ $ len$(walks) < max\_length \ \textbf{do}$
		\STATE \hspace{1cm} Select node $w$ uniformly at random from neighbors $\mathcal{N}(curr\_node)$ of $curr\_node$
		\STATE \hspace{1cm} Generate uniformly at random a number $0 \le p \le 1$
		\STATE \hspace{1cm} \textbf{if} $p < deg(curr\_node)/deg(w)$ \textbf{then}
		\STATE \hspace{1.5cm} $curr\_node = w$
		\STATE \hspace{1.5cm} $\textit{walks}$.append($curr\_node$)
		\STATE \hspace{1cm} \textbf{else}
		\STATE \hspace{1.5cm} stay at $curr\_node$
		\STATE \hspace{0.5cm} $\textbf{end}$ $\textbf{while}$
	\end{algorithmic}
	\label{alg:spatial_sample}
\end{algorithm}

\begin{figure}[h!t]
	\centering
	\includegraphics[width=0.8\linewidth]{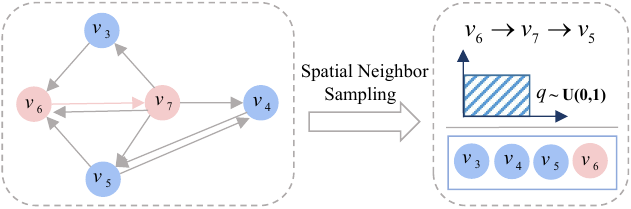}
	\caption{Schematic diagram of spatial space sampling based on Metropolis Hastings sampling. As shown in the figure, given node $v_{7} $, the spatial neighbor set $\{v_{3}, v_{4},v_{5},v_{6}\}$ is first selected. Then, select a node from the set of spatial neighbors according to the uniform distribution, calculate the acceptance rate using the proposal distribution and target distribution, and then sample to obtain the next node $v_{5} $.}
	\label{fig:spatial}
\end{figure} 

\subsection{Spatiotemporal Coupling Random Walks Encoding}
\textbf{Skip-Gram Embedding}
Word2vec is a natural language processing (NLP) model that is used for word embedding, which is used to represent words in a continuous vector space \cite{word2vec}. Word2vec is an important tool, since it has been widely used in many scenarios, such as machine translation, text classification, and sentiment analysis. Skip-gram is a type of the Word2vec, which learns to predict the context words that are most likely to appear around a given target word. 
Vector representations learned using the Skip-gram model have been found to capture various aspects of semantic and syntactic relationships between words. Building upon the obtained random walk sequences in both the temporal and spatial dimensions, we integrate the two to obtain a novel spatiotemporal random walk sequence $\mathcal{W}_{ts}$. 
The $\mathcal{W}_{ts}$ sequence consists of three parts that are concatenated together: first the target link, followed by the temporal-biased and spatial-biased sequences. 
This combined sequence allows for the exploration of both spatial and temporal information in a coherent manner, providing a more comprehensive understanding of the underlying network structure.
The mathematical description of Skip-gram used in dynamic network representation learning is as follows. Given a spatiotemporal random walk sequence $\mathcal{W}_{ts}$, it maximizes the co-occurrence probability among the nodes that appear within a window, 
\begin{equation}
	\underset{f}{max} \ log P(\mathcal{W}_{ts} = \{ v_{i-\omega},...,v_{i+\omega} \} \setminus v_{i} | f(v_{i})),
\end{equation}
where $f: \mathcal{V} \to \mathbb{R}^{D}$ is the embedding function, $\omega$ is the context window size for optimization. We assume conditional independence of the nodes, and the optimization problem can be 
formulated as:
\begin{equation}
	P(\mathcal{W}_{ts} | f(v_{i})) = 	\prod_{v_{j} \in \mathcal{W}_{ts}} P(v_{j} | f(v_{i})).
\end{equation}

\subsection{Edge Representation}
Using the node embedding obtained from optimizing the transmission of information in the time dimension and fully exploring the surrounding spatial structure through MH sampling, we generate edge representations via binary operations. The process of generating representations of edges is shown in Algorithm \ref{alg:edge_representation}. Specifically, we create the final edge representation $\emph{z}_{uv}$ using the embedding representations $\emph{z}_{u}$ and $\emph{z}_{v}$ of nodes $u$ and $v$. We experiment with five different binary operations:
\begin{itemize}
	\item{Concatenation: Concatenation of two vectors is the process of joining them end-to-end to form a longer vector.}
	\item{Mean: The mean of two vectors is their element-wise average.}
	\item{L1: The L1 metric quantifies the absolute difference between elements of the encoding vectors of two nodes.}
	\item{L2: L2 measures the sum of the squared differences between the corresponding elements of two encoding vectors of nodes.}
	\item{Hadamard: Hadamard multiplication of two node encoding vectors of equal dimension is used to generate an edge encoding vector with the same dimension.}
\end{itemize}

\begin{algorithm}[H]
	\caption{Edge Representation.}
	\begin{algorithmic}
		\STATE 
		\STATE{\bf Input:}
		Input networks $\mathcal{G = (V, E_{T}, T)}$, start\_edge $e^{t}_{v,u}=(v,u,t)$, the number of walks for each node $R$,  binary operations $D$=\{\textit{Concatenation, Mean, L1,  L2, Hadamard}\}.\\
		\STATE{\bf Output:} Edge Representation $\emph{z}_{uv}$.
		\STATE \hspace{0.5cm} Acquire $curr\_node = u$
		\STATE \hspace{0.5cm} \textbf{For} i = 1 : R  \textbf{do}
		\STATE \hspace{1cm} $\mathcal{W}_{st}^{i}(u) = [(v,u), \mathcal{W}_{t}^{i}(u),\mathcal{W}_{s}^{i}(u)]$
		\STATE \hspace{0.5cm}  $\mathcal{W}_{st}(u) = [\mathcal{W}_{st}^{1}(u), \mathcal{W}_{st}^{2}(u),...,\mathcal{W}_{st}^{R}(u)]$
		\STATE \hspace{0.5cm} Carry out Skip-gram with $\mathcal{W}_{st}(w), w \in \mathcal{V}$
		\STATE \hspace{0.5cm} Obtain node embedding $\emph{z}_{u}$, $\emph{z}_{v}$
		\STATE \hspace{0.5cm} Choose one binary operation $\odot$ from $D$
		\STATE \hspace{0.5cm} Edge embedding representation $\emph{z}_{uv}=(\emph{z}_{u} \odot \emph{z}_{v})$.
	\end{algorithmic}
	\label{alg:edge_representation}
\end{algorithm}

\section{Experiments}
\subsection{Experimental Setup}
\subsubsection{\textbf{Datasets}}
We evaluate the performance of COM using eight real-world dynamic networks: Contact, Hypertext, Enron, Radslaw, Fb-forum, Bitcoin-Alpha, Email-eu, and Wiki-elec. Table \ref{tab:datasets} provides a summary of the detailed information for each dataset. 
\begin{itemize}
	\item{\textbf{Contact }\cite{Contact}: This is a human contact network, which was collected through wireless devices worn by individuals attending an event. A link between two individuals is established if one person contacts another at a specific timestamp.} 
	\item{\textbf{Hypertext}\cite{Hypertext}: This network depicts the face-to-face interactions among attendees at the ACM Hypertext Conference 2009. Each node represents a conference visitor, while each edge denotes a face-to-face contact that lasts for 20 seconds.}
	\item{\textbf{Enron }\cite{Enron}: The Enron email network consists of emails sent between employees of Enron. Nodes in the network are individual employees and edges are individual emails.}
	\item{\textbf{Radslaw }\cite{Radslaw}: This network consists of email exchanges between employees of medium-sized manufacturing companies. Nodes in the network represent individual employees, while edges between nodes represent the exchange of emails.}
	\item{\textbf{Fb-forum }\cite{Fb-forum}: The dataset is derived from an online student forum at the university of California, similar to Facebook. The network consists of users as nodes and their interactions as links, spanning over a period of five months or more.}
	\item{\textbf{Bitcoin-Alpha }\cite{Bitcoin-Alpha}: This is who-trusts-whom network of people who trade using Bitcoin on a platform called Bitcoin Alpha. Nodes represent users, and edges represent the source's rating for the target, ranging from -10 to +10 in steps of 1.}
	\item{\textbf{Email-eu }\cite{Email-eu}: The network was created using email data collected from a large European research institution. Each node in the network represents a user, and edges between the nodes represent emails exchanged between two users.}
	\item{\textbf{Wiki-elec }\cite{Wiki-elec}: The dataset contains all administrator elections and vote history data based on the latest complete dump of Wikipedia page edit history.}
	
\end{itemize}

\begin{table}[]
	\centering
	\caption{Summary of Datasets Used in Our Experiments}
	\label{tab:datasets}
	\begin{tabular}{@{}lcccc@{}}
		\toprule
		Dataset & $|\mathcal{V}|$ & $|\mathcal{E_{T}}|$ & \textbf{$\bar{d}$} & Timespan(days) \\ \midrule
		Contact          & 274        & 28244        & 206.2      & 3.97                    \\
		Hypertext        & 113        & 20818        & 368.5      & 2.46                    \\
		Enron            & 151        & 50571        & 669.8      & 1137.55                 \\
		Radslaw          & 167        & 82927        & 993.1      & 271.19                  \\
		Email-eu         & 3783       & 24186        & 674.1      & 803.93                  \\
		Fb-forum         & 986        & 332334       & 75.0       & 164.49                  \\
		Bitcoin-Alpha    & 899        & 33720        & 12.8       & 1901.0                  \\
		Wiki-elec        & 7118       & 107071       & 30.1       & 1378.34                 \\ \bottomrule
	\end{tabular}
\end{table}

\subsubsection{\textbf{Evaluation Metrics}}
The accuracy of link prediction methods in continuous-time dynamic networks can be evaluated using various metrics. In this paper, we apply 1) AUC \cite{AUC} and 2) AP \cite{AP} , which measure the performance of the prediction method in terms of true positive rate, false positive rate, and precision.

\begin{itemize}
	\item{Area Under the Curve (AUC): AUC represents the accuracy of a machine learning algorithm from a general perspective, and is calculated as the area under the Receiver Operating Characteristic (ROC) curve. It measures the probability of a similarity value of a randomly chosen new link being greater than a randomly chosen nonexistent link. A higher AUC value indicates better model performance.}
	\item{Average Precision(AP): AP is a performance metric used for binary classification problems, which calculates the area under the precision-recall curve. It measures the quality of the model's positive predictions by considering both precision and recall, where precision is the fraction of true positives among the predicted positives, and recall is the fraction of true positives among all actual positives.}
	
\end{itemize}

\subsubsection{\textbf{Baselines}}
To evaluate the effectiveness of COM, we compare it with three kinds of methods: 1) heuristic methods including Common Neighbors (CN) and Jaccard Coefficient (JC), 2) static link prediction methods including node2vec, and LINE, and 3) dynamic link prediction methods including CTDNE, DySAT, and TGAT.

\begin{itemize}
	\item{\textbf{CN }\cite{CN}: The algorithm uses the number of common neighbors as an indicator to measure the possibility of establishing a link between two nodes.} 
	\begin{equation}
		CN(x, y) = |N(x) \cap N(y) |
	\end{equation}
	\item{\textbf{JC }\cite{JC}: This algorithm evaluates the probability of connecting edges also by measuring the number of common neighbors, it is the normalized version of CN.}
	\begin{equation}
		JC(x, y) = \frac{|N(x) \cap N(y) |}{|N(x) \cup N(y) |}
	\end{equation}
	\item{\textbf{Node2vec }\cite{Node2vec}: Node2vec is a representation learning algorithm for networks that learns low-dimensional representations of nodes by optimizing a second-order proximity measure using a biased random walk. Depth-first and breadth-first are controlled by parameters to learn the representation of nodes in the graph.}
	\item{\textbf{LINE }\cite{Line}: LINE is a scalable framework for learning high-quality node embeddings in large-scale networks. It optimizes a carefully designed objective function based on both first-order and second-order proximity.}
	\item{\textbf{CTDNE }\cite{CTDNE}: CTDNE is a method for learning low-dimensional representations of nodes in dynamic networks based on their temporal interactions and dynamics, using a node sampling strategy incorporating time information.}
	\item{\textbf{DySAT }\cite{Dysat}: DySAT is a dynamic graph embedding method that leverages self-attention to capture the temporal dependencies of nodes and their surrounding contexts for link prediction in dynamic networks.}
	\item{\textbf{TGAT }\cite{TGAT}: TGAT is a deep learning-based model for dynamic graph representation learning that uses attention mechanisms to capture temporal and structural information, achieving excellent performance on various dynamic graph prediction tasks.}
	
\end{itemize}


\begin{table*}[!ht]
	\centering
	\begin{center}
		\caption{AUC of different methods on different datasets}
		\label{Con_AUC}
	\begin{tabular}{llcccccccc}
		\toprule[0.8pt]
		Class   & Method   & Contact & Hypertext & Enron & Radslaw & \multicolumn{1}{l}{Email-eu} & \multicolumn{1}{l}{Fb-forum} & \multicolumn{1}{l}{Bitcoin-Alpha} & \multicolumn{1}{l}{Wiki-elec} \\ \hline
		\multicolumn{1}{c}{\multirow{2}{*}{\textbf{heuristic}}} 
		& CN \cite{CN}  & 0.969 & 0.761 & 0.936 & 0.942 & 0.962  &0.700 & 0.722 &0.780            \\
		\multicolumn{1}{c}{}   
		& JC\cite{JC}       & 0.901   & 0.788     & 0.883 & 0.862   & 0.965 & 0.684  & 0.792   & 0.844     \\ \hline
		\multirow{2}{*}{\textbf{static}}                      
		& Node2vec \cite{Node2vec} & 0.878   & 0.693  & 0.721 & 0.794   & 0.717     & 0.722      & 0.836   & 0.601             \\
		& Line \cite{Line}    & 0.736   & 0.621     & 0.550  & 0.615   & 0.650      & 0.648           & 0.672       & 0.625     \\ \hline
		\multirow{4}{*}{\textbf{dynamic}}                      
		& CTDNE \cite{CTDNE}    & 0.907   & 0.893     & 0.843 & 0.841 & 0.729   & 0.818     & 0.837   & 0.826                           \\
		& DySAT \cite{Dysat}    & 0.947   & 0.715     & 0.878 & 0.817   & 0.959                        & 0.821                        & 0.846                             & 0.816                         \\
		& TGAT \cite{TGAT}    & 0.921   & \textbf{0.959}     & 0.786 & 0.905   & 0.719                        & 0.878                        & 0.872                             & 0.953                         \\
		& COM  & \textbf{0.989}   & \textbf{0.959}     & \textbf{0.953} & \textbf{0.986}  & \textbf{0.979}  &\textbf{0.923}                        & \textbf{0.878}   & \textbf{0.956}                         \\ \toprule[0.8pt]
	\end{tabular}
	\end{center}
\end{table*}

\begin{table*}[!ht]
	\begin{center}
		\caption{AP of different methods on different datasets}
		\label{Con_AP}
	\begin{tabular}{llcccccccc}
		\toprule[0.8pt]
		Class & Method  & Contact & Hypertext & Enron & Radslaw & \multicolumn{1}{l}{Email-eu} & \multicolumn{1}{l}{Fb-forum} & \multicolumn{1}{l}{Bitcoin-Alpha} & \multicolumn{1}{l}{Wiki-elec} \\ \hline \multicolumn{1}{c}{\multirow{2}{*}{\textbf{heuristic}}} 
		& CN\cite{CN}       & 0.973          & 0.753          & 0.933          & 0.935          & 0.959                         & 0.687                        & 0.721                             & 0.776                         \\
		\multicolumn{1}{c}{}                                    
		& JC\cite{JC}       & 0.817          & 0.754          & 0.867          & 0.820           & 0.964                        & 0.616                        & 0.748                             & 0.825                         \\ \hline
		\multirow{2}{*}{\textbf{static}}                        
		& Node2vec\cite{Node2vec} & 0.835          & 0.662          & 0.703         & 0.772          & 0.709                        & 0.735                        & 0.847                             & 0.737                         \\
		& Line \cite{Line}    & 0.718          & 0.607          & 0.536          & 0.584          & 0.649                        & 0.612                        & 0.724                             & 0.619                         \\ \hline
		\multirow{4}{*}{\textbf{dynamic}}                       
		& CTDNE \cite{CTDNE}   & 0.918           & 0.897           & 0.837        & 0.846          & 0.691                        & 0.825                       & 0.820                             & 0.726                         \\
		& DySAT\cite{Dysat}    & 0.918          & 0.709          & 0.855          & 0.798          & 0.951                        & 0.803                        & 0.822                             & 0.804                         \\
		& TGAT\cite{TGAT}     & 0.906          & \textbf{0.927} & 0.774          & 0.882          & 0.691                        & 0.831                        & 0.676                             & \textbf{0.956}                \\
		& COM  & \textbf{0.982} & 0.923          & \textbf{0.946} & \textbf{0.982} & \textbf{0.976}               & \textbf{0.893}               & \textbf{0.863}                    & 0.938                         \\ \toprule[0.8pt]
	\end{tabular}
	\end{center}
\end{table*}
We assess the effectiveness of the COM on the temporal link prediction task. For this, we sort edges in each graph by ascending chronological order and use the first 75\% for representation learning. We consider the remaining 25\% as positive links and randomly sample an equal number of negative edges for generating labeled examples for link prediction. Node2vec, LINE, and DySAY use the logistic regression classifier; CTDNE and TGAT use the MLP classifier. We follow the setting in the original paper for the choice of parameters in the baseline experiment. The specific settings of parameters in the COM model are as follows, the window size $\omega$ in Skip-gram model is set to 10. Besides, the maximum walk length $L$ is set to 80, the number of walks for each node $R$ is set to 10, and the embedding dimension $D$ is set to 128. 

\subsection{Experimental Results}
The link prediction performance of all methods was evaluated in terms of AUC and AP, and the results are summarized in Table \ref{Con_AUC} and Table \ref{Con_AP}. The reported performance is an average over ten repetitions. Overall, our proposed model shows competitive performance compared to other methods.

The table shows that the heuristic method performs well on dense networks, while static models such as Node2vec and LINE perform poorly since they only consider node representations obtained from random walks without time information. Dynamic graph models, which capture time information in the network, outperform static models. The performance of CTDNE is improved by using a random walk algorithm that considers time information. On most datasets, graph neural network methods DySAT and TGAT perform better than other methods, as they can aggregate neighbor information through attention mechanisms to obtain better node representations. Our proposed method, which optimizes temporal information propagation while coupling spatial structures, achieves better performance in most cases than other dynamic models. Although its AP values perform slightly worse in Hypertext and Wiki-elec, it is still competitive. In summary, our method can achieve satisfactory results on both sparse and dense graphs. The promising results presented above suggest that the node embeddings obtained through our method can be utilized as beneficial features for link prediction tasks.
\begin{figure*}
	\centering
	\includegraphics[width=0.95\linewidth]{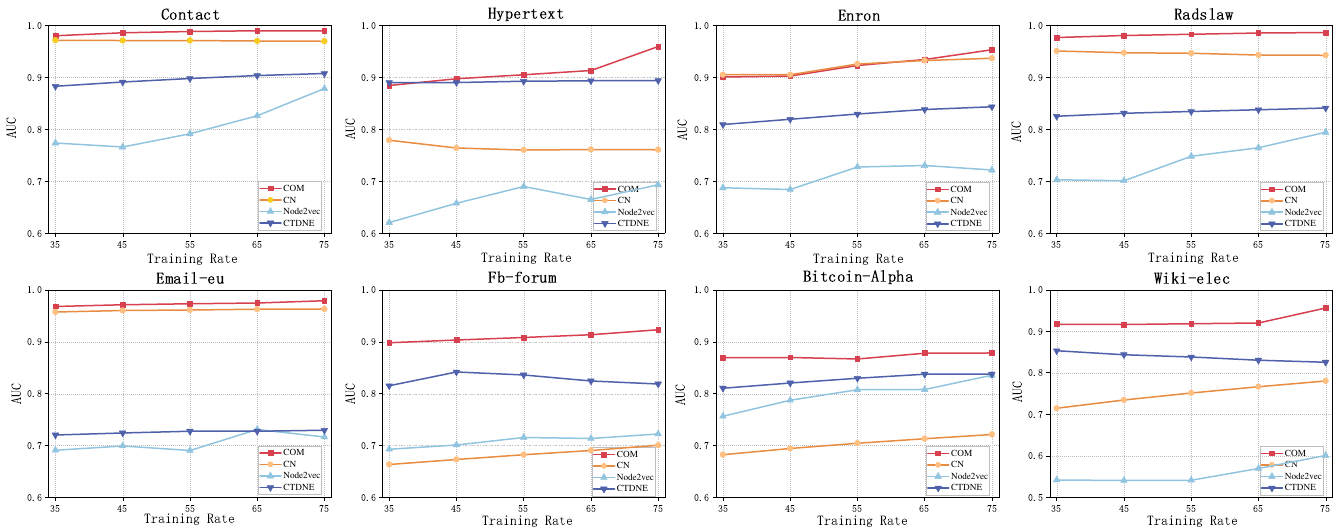}
	\caption{AUC comparison was performed on all datasets of CN, Node2vec, CTDNE and COM using different percentage training links. For each dataset, 35\%, 45\%, 55\%, 65\%, and 75\% of all links in networks were used as training sets.}
	\label{fig:training_rate}
\end{figure*}

\begin{figure}
	\centering
	\includegraphics[width=0.9\linewidth]{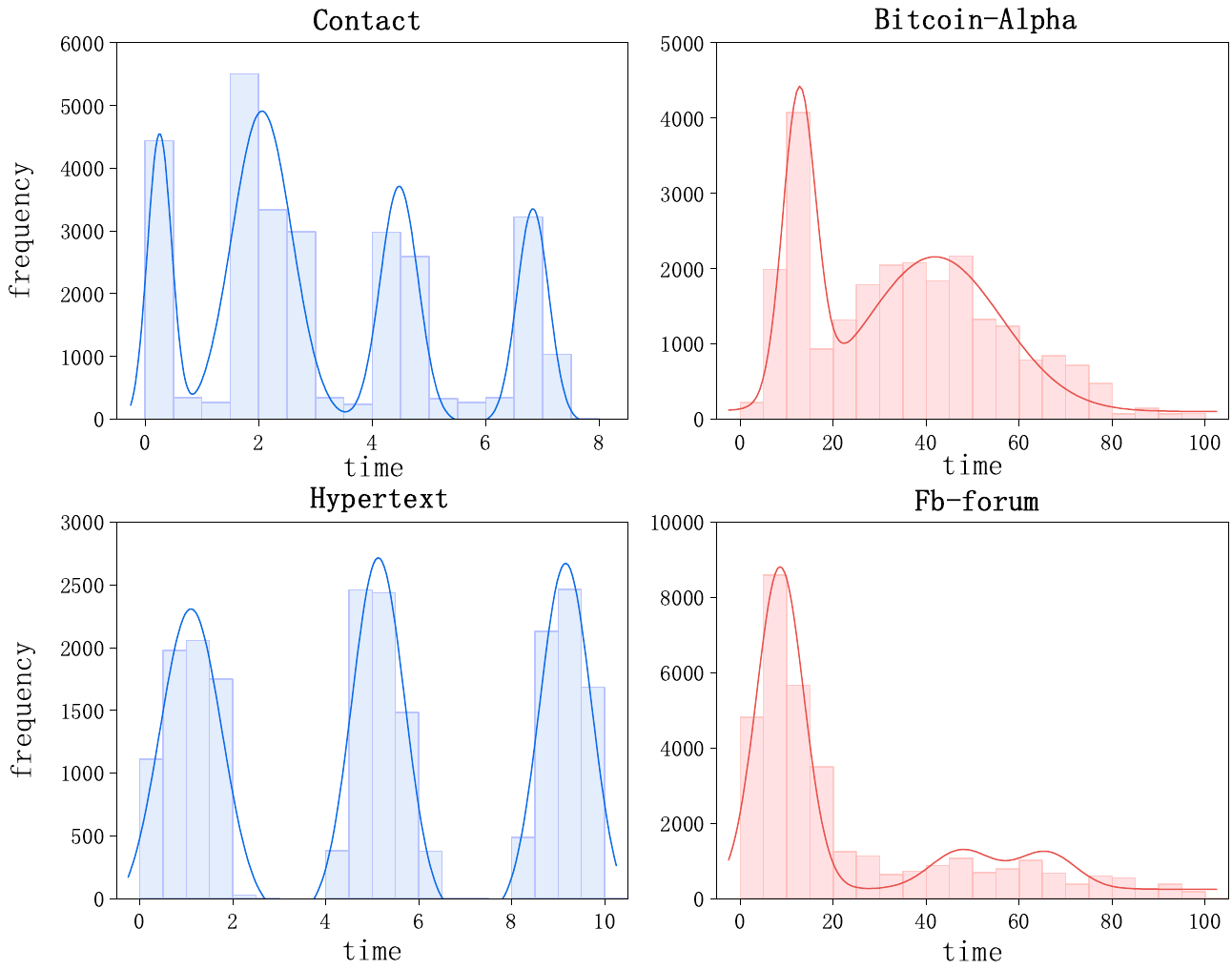}
	\caption{Network growth pattern of different datasets. We normalize the time of the network to start from 0. The horizontal axis represents the evolution time of the network, and the vertical axis represents the frequency of edges appearing during the corresponding time period.}
	\label{fig:frequency}
\end{figure}
We experiment with different training percentages and presented the AUC results of our model comparing CN, node2vec, and CTDNE in Fig. \ref{fig:training_rate}.
We used the last 25\% of the edges as the test set and gradually increased the number of training edges in reverse chronological order from 35\% to 75\%. Our proposed method's AUC values are indicated by red lines, and baseline methods are denoted with different colors. The figure showed that our proposed method outperformed almost all baseline methods under different percentages of training data.  The CN method's performance deteriorated with the increase of the training set proportion on the Contact, Hypertext, and Radslaw datasets, possibly due to the addition of edges far from the current time period to the network, which introduced more noise to the model. Node2vec's performance is unstable and fluctuates significantly across most datasets. The trend of CTDNE is very similar to our model.
In general, a larger training set should generally provide a better performance. 
COM exhibits consistent and stable performance on different training sets of most datasets, demonstrating its capacity to learn essential network features and patterns and its robustness and generalization across different data scales.
Notably, for datasets Hypertext, Enron, and Wiki-elec, when increasing the training rate from 65\% to 75\% , there is a notable enhancement in the model's performance, which emphasizes the relevance of recent data for precise predictions.
\begin{figure}[h!b]
	\centering
	\includegraphics[width=0.8\linewidth]{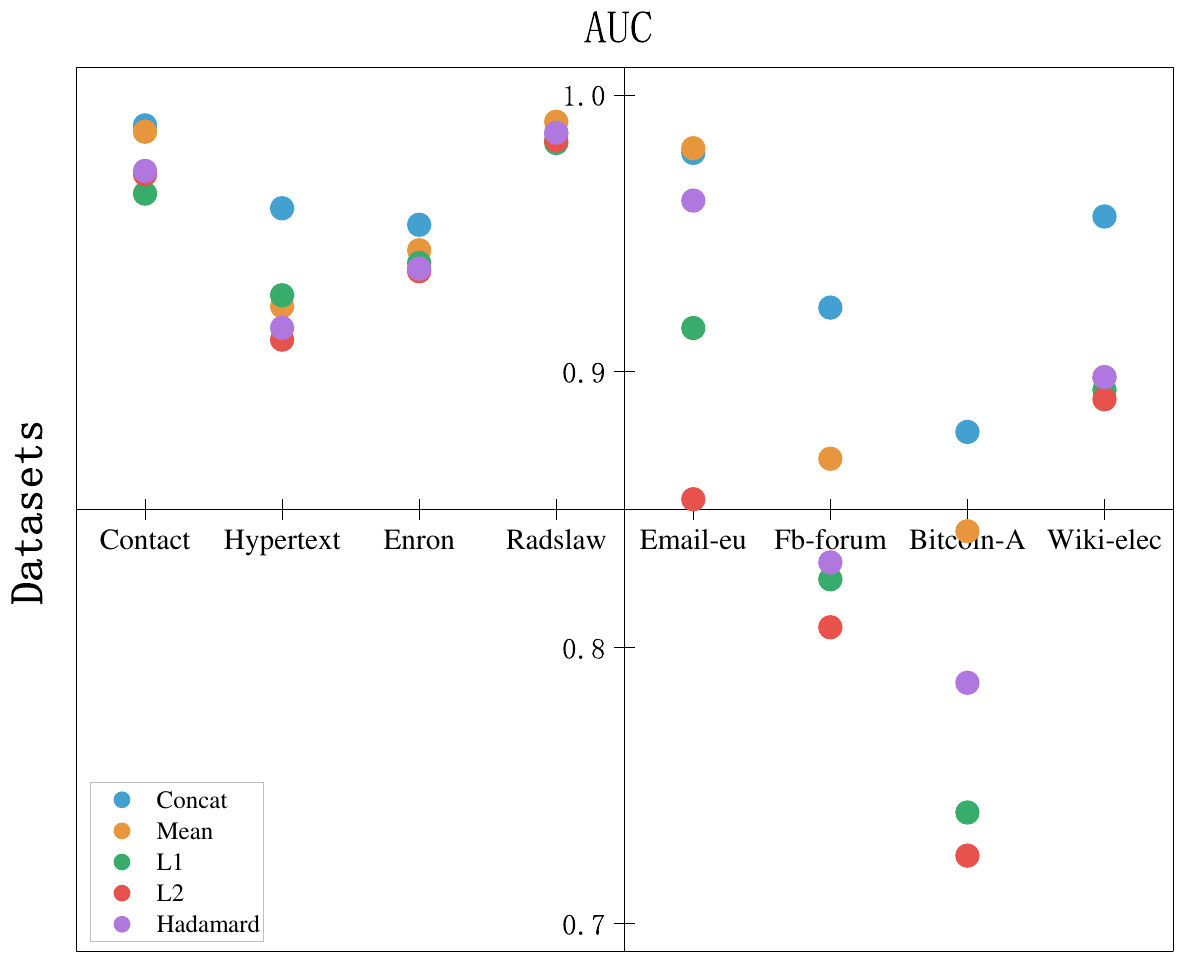}
	\caption{Different edge representation methods test in COM across all datasets, measured by AUC.}
	\label{fig:edge-representation}
\end{figure}

Additionally, we conducted a comprehensive analysis of the model's performance across different datasets. To provide further insights into these differences, we visualized the data distribution of edges in various networks, as shown in Fig. \ref{fig:frequency}, and fitted them with curves. Notably, networks Contact and Hypertext, depicted in blue, showcased superior performance on our proposed COM model, while networks Bitcoin-Alpha and Fb-forum, represented in red, exhibited relatively poor performance. The figure reveals that the evolution and growth pattern of the network that performs well in the model shows regularity.
Conversely, networks with lower performance display erratic distribution and higher volatility, with noticeably slower growth rates in later stages. These observations further confirm the effectiveness of the COM model in capturing information pertaining to the growth and evolution patterns of dynamic networks.

\subsection{Embedding Analysis}
We use binary operations to combine node representations into edge representations and compare the performance of different methods. Fig. \ref{fig:edge-representation} shows that the effectiveness of each method varies across datasets. The experimental results show that the concatenation operation outperforms other edge representation methods on most datasets, except for Radslaw and Eamil-eu where the Mean operation dominates. The Mean operation is  inferior, while L2 performs poorly across most datasets. Based on these results, it can be concluded that other binary operations lead to information loss when compared to concatenation. Thus, we selected concatenation as the binary operation for our model. Notably, the performance of different methods on the Forum and Bitcoin-Alpha datasets shows significant variability, possibly due to the network's sparsity and the slow growth of the network in the later stage as shown in Fig. \ref{fig:frequency}.

\begin{figure*}[h!t]
	\centering
	\includegraphics[width=\linewidth]{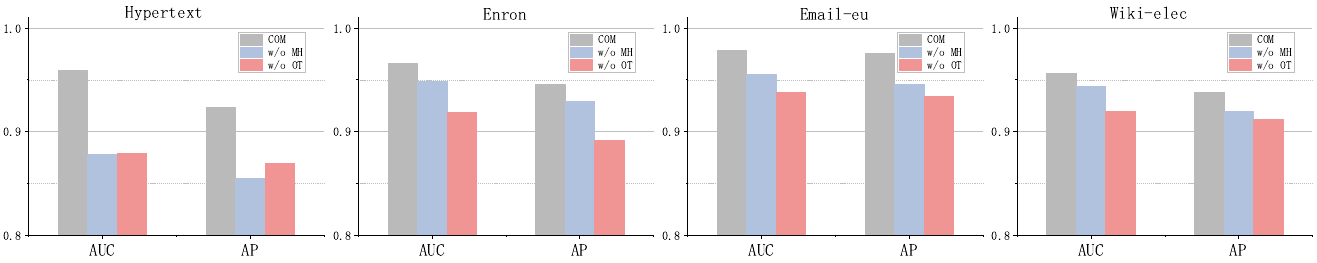}
	\caption{Comparative analysis of AUC and AP Results. Evaluating the impact of only considers optimal temporal information transmission (\textit{w/o OT}) and spatial structure information (\textit{w/o MH}) on COM performance, respectively.}
	\label{fig:ablation}
\end{figure*}

\subsection{Ablation Study}
In Fig. \ref{fig:ablation}, we present the results of ablation studies conducted on the Hypertext, Enron, Email-eu and Wiki-elec in both temporal and spatial dimensions to examine the impact of temporal-biased and spatial-biased sampling strategy on dynamic link prediction. Specifically, we additionally demonstrate the ablation test in terms of discarding the temporal validity and shortest distance of information transmission based on optimal transmission when extracting spatial structure  information (denoted as \textit{w/o OT}), and removing the extraction of spatial structure based on Metropolis Hastings when sampling temporal-biased walks (denoted as \textit{w/o MH}).
Our findings indicate that considering only temporal or spatial information alone leads to a decrease in link prediction performance. This highlights the effectiveness of our spatiotemporal biased walking sampler, and provides valuable insights into the intricate interplay between temporal dynamics and spatial structure in complex networks. Furthermore, when the optimal transmission of information between nodes is not considered in sampling, the prediction performance decreases significantly, demonstrating the importance of our proposed temporal-biased strategy in reducing the loss of information propagation in the network.

\section{Conclusion}
In this paper, we introduce a dynamic link prediction framework based on random walk, which addresses the issue of information loss in time dimension sampling by utilizing the optimal transmission algorithm. To efficiently integrate spatiotemporal information, we utilize the MH sampling algorithm to extract high-order structures around the target link. Experiments on eight datasets from different fields show that COM outperforms a wide range of methods in both heuristic and graph neural networks methods.
In future work, we aim to enhance our model's performance by integrating external information, such as text data and user attributes. Moreover, developing a multi-task link prediction model that can predict different types of links simultaneously will provide more comprehensive insights into network dynamics and structure.

\section*{Acknowledgments}
This work was supported by the National Natural Science Foundation of China (Grant Nos. 62276013, 62141605, 62050132), the Beijing Natural Science Foundation (Grant No. 1192012), and the Fundamental Research Funds for the Central Universities.

\bibliographystyle{IEEEtran}
\bibliography{refe}

\section{Biography Section}
\begin{IEEEbiography}[{\includegraphics[width=1in,height=1.25in,clip,keepaspectratio]{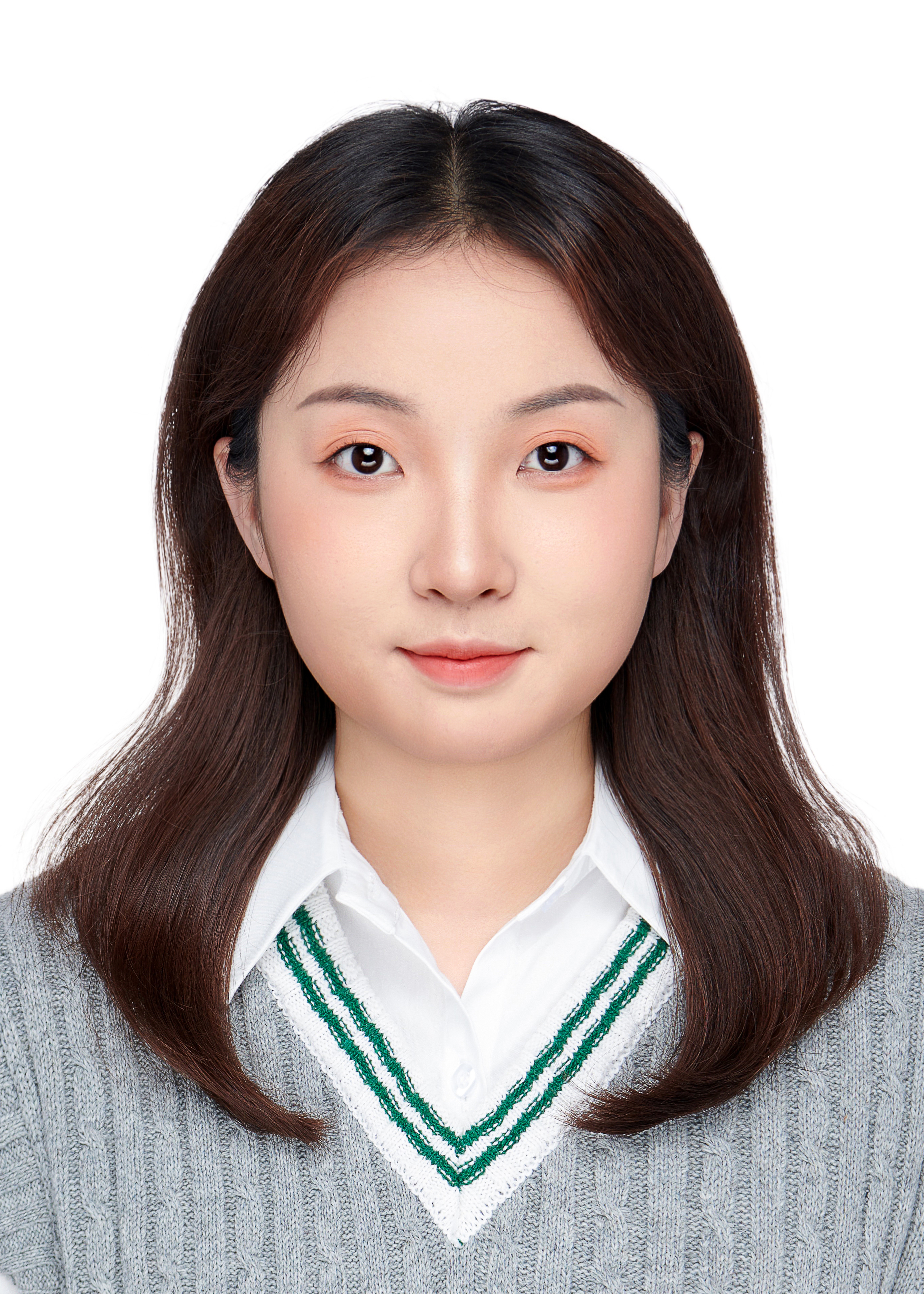}}]{Ruizhi Zhang}
	is currently pursuing the Ph.D. degree with the School of Mathematical Sciences, Beihang University, Beijing, China. Her current research interests include complex networks, link prediction and graph representation learning.
\end{IEEEbiography}

\begin{IEEEbiography}[{\includegraphics[width=1in,height=1.25in,clip,keepaspectratio]{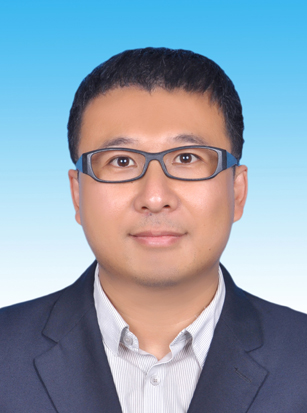}}]{Wei Wei}
	received the Ph.D. degree in mathematics from the School of Mathematical Sciences, Peking University, Beijing, China, in 2009. He is currently an Associate Professor with the School of Mathematical Sciences, Beihang University, Beijing, China. His research interests include dynamical system and complexity, complex networks, and artificial intelligence.
\end{IEEEbiography}

\begin{IEEEbiography}[{\includegraphics[width=1in,height=1.25in,clip,keepaspectratio]{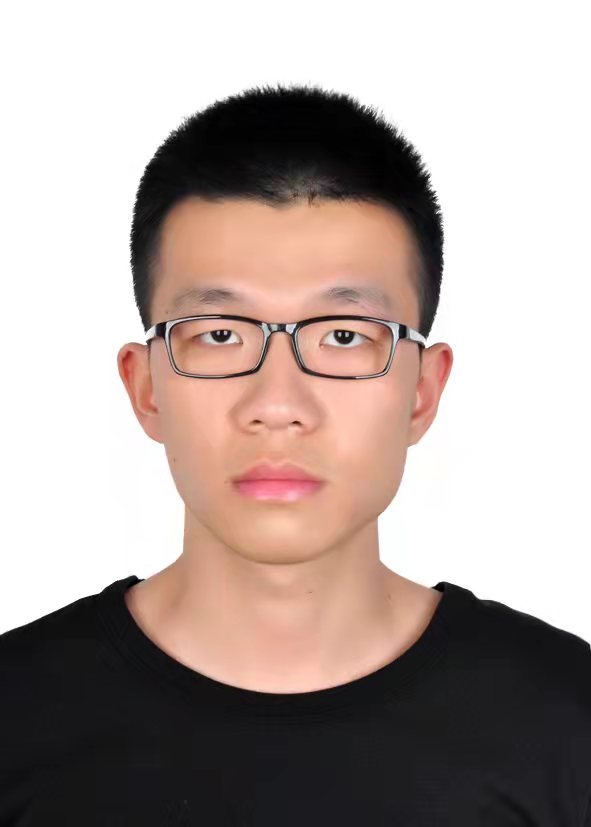}}]{Qiming Yang}
	received the B.Sc. degree from Beihang University, Beijing, China, in 2021. He is currently pursuing the master's degree with the School of Mathematical Sciences, Beihang University, Beijing, China. His current research interests include complex networks, link prediction and graph representation learning.
\end{IEEEbiography}

\begin{IEEEbiography}[{\includegraphics[width=1in,height=1.25in,clip,keepaspectratio]{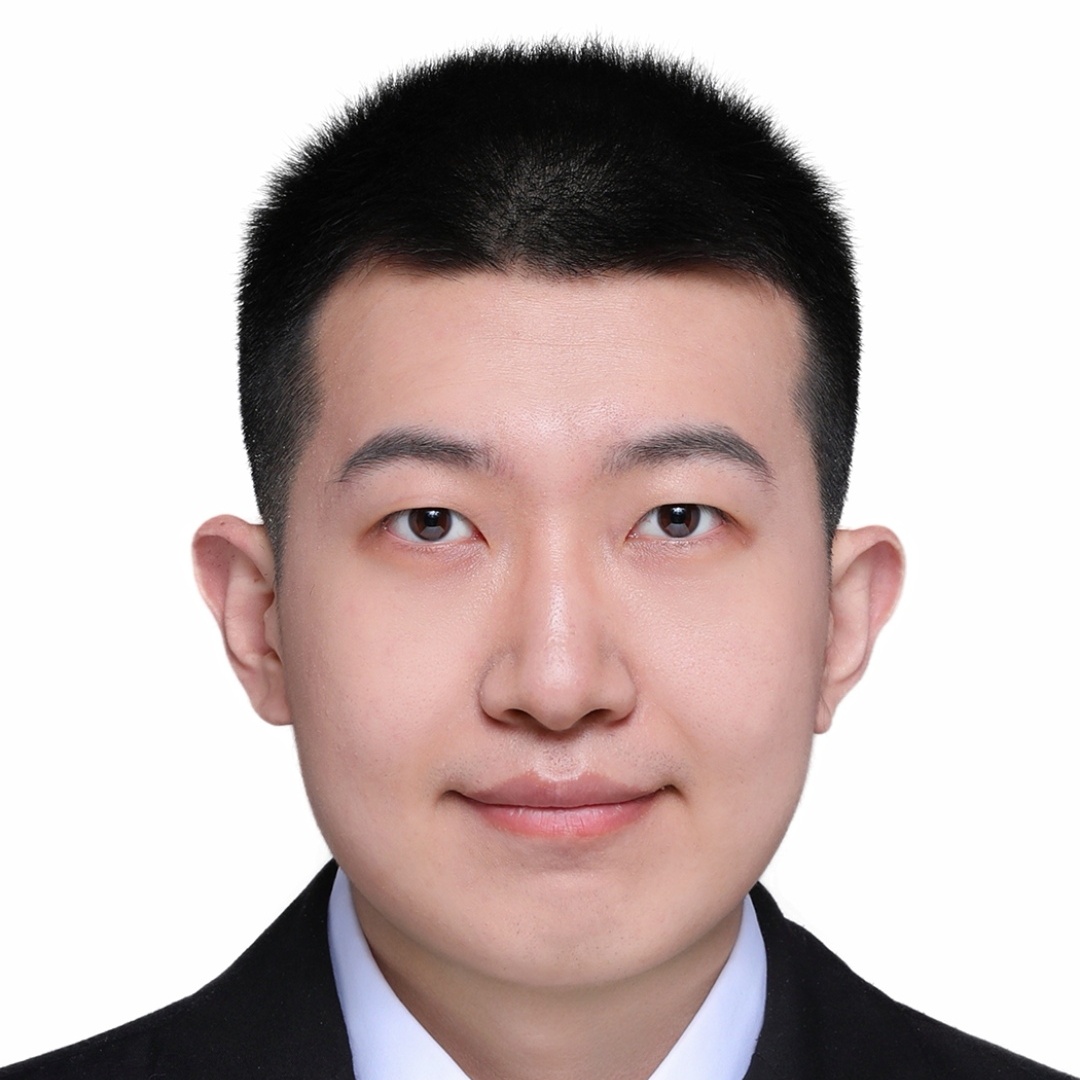}}]{Zhenyu Shi}
	received the B.Sc. degree from Beihang University, Beijing, China, in 2019. He is currently pursuing the Ph.D.’s degree with the School of Mathematical Sciences, Beihang University, Beijing, China and studying in Kyushu University as an exchange student.  His current research interests include complex networks and evolutionary game theory. 
\end{IEEEbiography}

\begin{IEEEbiography}[{\includegraphics[width=1in,height=1.25in,clip,keepaspectratio]{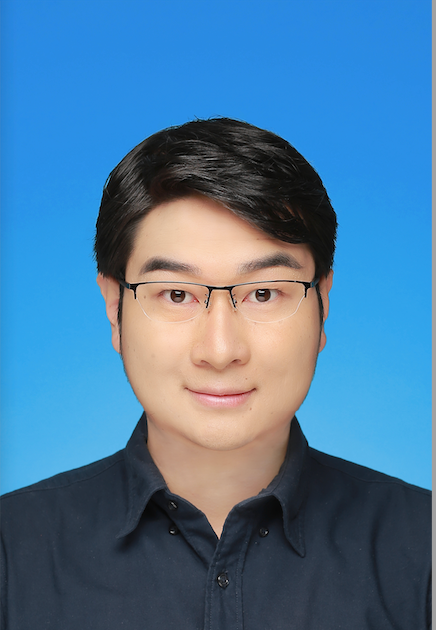}}]{Xiangnan Feng}
	received the Ph.D. degree in mathematics from the School of Mathematical Sciences, Beihang University, Beijing, China, in 2020. He is currently a Postdoctoral Fellow with the Complexity Science Hub, Vienna, Austria. His research interests include complex networks, computing social science, and artificial intelligence.
\end{IEEEbiography}

\begin{IEEEbiography}[{\includegraphics[width=1in,height=1.25in,clip,keepaspectratio]{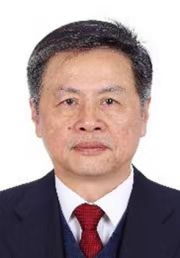}}]{Zhiming Zheng} received his Ph.D. degree in mathematics from the School of Mathematical Sciences, Peking University, Beijing, China, in 1987. Currently, he is a Professor with the Institute of Artificial Intelligence, Beihang University, Beijing, China. He is the Member of Chinese Academy of Sciences. His research interests include refined intelligence theory, information security, and complex information systems.
\end{IEEEbiography}

\end{document}